\documentclass[12pt,a4paper]{article}
\usepackage{amsmath}

\usepackage{amssymb,latexsym}

\newtheorem{theorem}{Theorem} [section]

\newtheorem{corollary}{Corollary} [section]

\title{Eigenvalues of the Laplacian acting on $p$-forms and metric conformal deformations }
\author{Bruno Colbois and Ahmad El Soufi}
\date{}

\begin{document}

\maketitle

\bigskip
\maketitle

\bigskip

\begin{abstract}

Let $(M,g)$ be a compact connected orientable Riemannian manifold of dimension $n\ge4$ and let
$\lambda_{k,p} (g)$ be the $k$-th positive eigenvalue of the Laplacian $\Delta_{g,p}=dd^*+d^*d$ acting on differential forms of degree $p$ on $M$. 
We prove that the metric $g$ can be conformally deformed to a metric $g'$, having the same volume as $g$, with arbitrarily large $\lambda_{1,p} (g')$ for all $p\in[2,n-2]$.

Note that for the other values of $p$, that is $p=0, 1, n-1$ and $n$, one can deduce from the literature that, $\forall k >0$, the $k$-th eigenvalue $\lambda_{k,p}$ is uniformly bounded on any conformal class of metrics of fixed volume on $M$.

For $p=1$, we show that, for any positive integer $N$, there exists a metric $g_N$ conformal to $g$ such that, $\forall k\le N$, $\lambda_{k,1} (g_N) =\lambda_{k,0} (g_N)  $, that is, the first $N$ eigenforms of $\Delta_{g_N,1}$ are all exact forms.

\end{abstract}

\bigskip


\medskip


\bigskip
\bigskip
\bigskip
\noindent
B. Colbois~: Universit\'e de Neuch\^atel, Laboratoire de
Math\'ematiques, 13 rue E. Argand,
2007 Neuch\^atel, Switzerland.

\smallskip
\noindent
E-mail: Bruno.Colbois@unine.ch

\bigskip
\noindent
A. El Soufi~: Universit\'e de Tours, Laboratoire de Math\'ematiques
et Physique Th\'eorique, UMR-CNRS 6083, Parc de Grandmont, 37200 Tours, France.

\smallskip
\noindent
E-mail: elsoufi@univ-tours.fr.
\eject

\section{Introduction and Statement of Results}\label{1}

Let $M$ be a closed connected orientable differentiable manifold of dimension $n \geq 2$. For any Riemannian metric $g$ on $M$, we denote by $\Delta_{g,p}=dd^*+d^*d$ the corresponding Laplacian acting on differential forms of degree $p$ on $M$. The spectrum $spec(\Delta_{g,p})$ of $\Delta_{g,p}$ consists of an unbounded sequence of nonnegative eigenvalues which starts from zero if an only if the $p$-th Betti number $b_p(M)$ of $M$ does not vanish. The sequence of positive eigenvalues of $\Delta_{g,p}$ is denoted by
$$ 0 < \lambda_{1,p}(g) \leq \lambda_{2,p}(g)
\cdots  \leq \lambda_{k,p} (g) \leq \cdots \to \infty,$$ 
The case $p=0$ corresponds to the Laplace-Beltrami operator acting on functions.

\medskip

The study of the behavior of the eigenvalues of the Laplacian under metric deformations constitutes one of the main topics in spectral geometry. In particular, it is of great interest to find situations where the behavior of $\lambda_{k,p}(g) $ for $p\ge1$, differs completely from that of $\lambda_{k,0}(g) $ (see for instance the situation 
of collapsing manifolds described in [CC] and [L]). 
 
\medskip

 In what follows, $\lambda_{k,p}$ will be considered as a functional on the space ${\cal M}(M)$ of all Riemannian metrics of {\it
volume one} on $M$ (the "volume one" condition is a normalization since $\lambda_{k,p}(cg)={1\over c}\lambda_{k,p}(g)$). We will also consider its restriction to a given conformal class of metrics of fixed volume $[g]= \left\{ g'\in {\cal M}(M) \ |\  g'\  {\hbox {is conformal to}} 
\ g \right\}$.
 
\medskip

There exists an extensive literature dealing with the scalar case $p=0$ (see for instance the introduction of [CE] for a synthetic presentation of the main results). In particular, it is well known that, 
\begin{itemize}

\item[(i)] if $n=2$, then $\lambda_{k,0}$ is bounded on ${\cal M}(M)$. More precisely, if $M$ is an orientable compact surface of genus $\gamma$, then,  there exists a universal constant $C$ such that $\forall k\ge1$,  
\begin{equation}
\sup _{g\in{\cal M}(M)} \lambda_{k,0}(g) \leq C(\gamma + 1)k
\end{equation}
(see [YY] for $k=1$ and [K] for general $k$). 

 \item[(ii)] if $n \ge3$, then $\lambda_{1,0}$ is not bounded on ${\cal M}(M)$, that is (see [CD])
 \begin{equation}
\sup _{g\in{\cal M}(M)} \lambda_{1,0}(g)=\infty.
\end{equation}

\item[(iii)] $\forall g \in {\cal M}(M)$, $\lambda_{k,0}$ is bounded on $[g]$. Indeed, there exists a constant  $C([g])$ such that, $\forall k\ge1$,  
\begin{equation}
\sup _{g'\in[g]} \lambda_{k,0}(g') \leq C([g]) k^{2/n}
\end{equation} 
 (see [EI] for $k=1$ and [K] for general $k$). 
\end{itemize}

\medskip

\ {\it Problem : do these results still hold when we replace $ \lambda_{k,0}$ with $ \lambda_{k,p}$, $p\ge1$?}

\medskip

Using the Hodge-de Rham decomposition theorem, one can easily see that for an orientable compact surface, $\lambda_{k,0}(g)=\lambda_{k,1}(g)=\lambda_{k,2}(g)$. Thus, the result (i) above is also valid for the eigenvalues of the Laplacian on forms.

\medskip

The case $n\ge3$ and $p\ge1$ was first considered by Tanno [T] who examined a 1-parameter family of Berger metrics $\left(g_t\right)_t$ on the sphere ${\mathbb S}^3$ for which $ \lambda_{1,0}(g_t)\to \infty$ as $t\to \infty$, and showed that $ \lambda_{1,1}(g_t)$ is bounded. On the other hand, Gentile and Pagliara [GP] proved that if $n \ge4$, then for any $p\in [2,n-2]$, $\lambda_{1,p}$ is not bounded on ${\cal M}(M)$ :
\begin{equation}
\sup _{g\in{\cal M}(M)} \lambda_{1,p}(g)=\infty.
\end{equation}
Hence, the result (ii) above is still valid for $\lambda_{1,p}$, $2 \le p \le n-2$. Note that the case $p=1$ (and its dual $p=n-1$) remains open.

\medskip

In this paper, we focus, as in the result (iii) above, on the restriction of the functional $\lambda_{k,p}$ to a given conformal class of metrics of fixed volume $[g]$ on $M$. As we saw in (iii) above, for any $k\ge1$, $\lambda_{k,0}$ is bounded on $[g]$. Since $d\Delta_{g,0 }=\Delta_{g,1}d$, the differential $df$ of any eigenfunction $f$ is an eigenform of degree 1 with the same eigenvalue. In particular, we have, $\forall k\ge 1$,
$$ \lambda_{k,1}(g)\le \lambda_{k,0}(g).$$
Consequently, $\lambda_{k,1}$ is also bounded on $[g]$ and the result (iii) is still valid for the eigenvalues of the Laplacian acting on 1-forms. Thanks to the Hodge duality,  the same is also true for the eigenvalues of the Laplacian on $(n-1)$-forms and on $n$-forms.
 Thus, {\it $\lambda_{k,p}$ is bounded on $[g]$ for $p=0$, $1$, $n-1$ and $n$}.

\medskip
Our main theorem shows that the situation differs drastically for the other values of $p$. Indeed, it turns out that if $2 \le p \le n-2$, then the functional $\lambda_{1,p}$ is never bounded on $[g]$, whatever the 
conformal class $[g]$ we consider. More precisely, one can conformally deform any metric $g$ on $M$ in such a way that the $n-3$ eigenvalues $\lambda_{1,p}$, $2 \le p \le n-2$, become simultaneously arbitrarily large. 

\begin{theorem} Let $(M,g)$ be a compact Riemannian manifold of dimension 
$n\ge 4$ and let $A$ be any positive real number. There exists a metric $g_A$ conformal to $g$ and having the same volume as $g$, such that, $ \forall p \in [2, n-2]$, 
$$\lambda_{1,p}(g_A)\ge A.$$
In particular, $ \forall p \in [2, n-2]$, we have
$$\sup _{g'\in [g]} \lambda_{1,p}(g')=\infty.$$
  \end{theorem} 
Or, equivalently,
$$\sup\{\lambda_{1,p}(g')V(g')^{2/n} | \ g' \ {\hbox {conformal to}} \ g\} = \infty,$$ 
where $V(g')$ denotes the volume of the manifold $(M,g')$. 

\medskip

Note that when the supremum is taken over all the metrics of volume 1, then there is no difference between the scalar case $p=0$ and the case $2 \le p \le n-2$ (that is, $\sup _{g\in{\cal M}(M)} \lambda_{1,p}(g)=\infty$) as proved in [GP]. 

\medskip
As an application of Theorem 1.1 and the result of Korevaar (3) above, we obtain the following. 

\begin{corollary} 
Let $(M,g)$ be a compact Riemannian manifold of dimension 
$n\ge 4$. For any positive integer $k$,  there exists a metric $g_k$ conformal to $g$ such that, $\forall$ $2 \le p \le n-2$, 
$$\lambda_{1,p} (g_k) > \lambda_{k,0} (g_k)$$
and
$$\lambda_{1,p} (g_k) > \lambda_{k,1} (g_k).$$
  \end{corollary} 

If we set, like in [Tk],  ${\hbox{Gap}}_{k,l}^{p,q}(g)=\lambda_{k,p}(g)-\lambda_{l,q}(g)$, then Corollary 1.1 tells us that, $\forall k\ge1$ and $ l\ge1$, every conformal class of metrics $[g]$ contains metrics $g'$ satisfying, $\forall$ $ 2 \le p \le n-2$, 
$${\hbox{Gap}}_{k,l}^{p,0}(g') >0\ \  {\hbox{and}}\  \ {\hbox{Gap}}_{k,l}^{p,1}(g')>0.$$

Recall that the gaps $ {\hbox{Gap}}_{k,k}^{0,1}(g)=\lambda_{k,0}(g)-\lambda_{k,1}(g)$ between the eigenvalues on functions and those on 1-forms are all nonnegative. An interesting question is to know whether the nature (zero or nonzero) of some of them is imposed by the differential or the conformal structures of the manifold. 

\medskip

Takahashi [Tk] (see also [GS]) proved that any manifold $M$ of dimension $n\ge3$ admits both a metric $g_1$ with ${\hbox{Gap}}_{1,1}^{0,1}(g_1)=0$ and a metric $g_2$ with ${\hbox{Gap}}_{1,1}^{0,1}(g_2)> 0$. Thus, the nature of the first gap is not an invariant of the differential structure of $M$.

\medskip

The following result shows that there is no differential or conformal obstruction for any finite number of gaps to be zero.
 Indeed,  we have the following. 

\begin{corollary} 
Let $(M,g)$ be a compact Riemannian manifold of dimension 
$n\ge 4$. For any positive integer $N$,  there exists a metric $g_N$ conformal to $g$ such that, $\forall k\le N$, 
$$\lambda_{k,1} (g_N) =\lambda_{k,0} (g_N),  $$
that is, the first $N$ eigenforms of $\Delta _ {g_N,1}$ are all exact forms.
  \end{corollary} 

This means that, $\forall N\in {\mathbb N} $, any metric $g$ can be conformally deformed to a metric $g_N\in [g]$ whose $N$ first gaps vanish, that is, 
$${\hbox{Gap}}_{1,1}^{0,1}(g_N)=\cdots = {\hbox{Gap}}_{N,N}^{0,1}(g_N) = 0.$$

\medskip
Finally, we point out that the existence of upper bounds on $ \lambda_{1,0}$ under other types of restrictions has been investigated, for instance in [BLY] and [P]. It would be of course interesting to study the existence of such upper bounds for $ \lambda_{1,p}$ with $p\ge 1$.

\section {Proof of Results} 

The idea of the proof  of Theorem 1.1 consists in  showing that the 
construction of [GP] to produce large eigenvalues for forms may be 
realized by conformal deformations. As in [CE], we will use the fact 
that any manifold is locally almost Euclidean. This implies that we 
can suppose first that the metric is Euclidean on a small 
neighbourhood of a point, construct a conformal deformation of the 
Euclidean metric, and then show that, by quasi-isometry, it gives a 
conformal deformation of the initial metric with the property we want. 

Before going into the proof, let us recall some basic facts about the Laplacian on forms. Indeed,
the Hodge-de Rham theorem gives the following orthogonal decomposition of the space $\Lambda^p(M)$ of differential forms of degree $p$ on M:
$$\Lambda^p(M)= Ker\Delta_{g,p}\oplus d\Lambda^{p-1}(M) \oplus d^*\Lambda^{p+1}(M).$$
The operator $\Delta_{g,p}$ leaves this decomposition invariant. Hence, if we denote by $\left\{\lambda'_{k,p}(g); k\ge1 \right\}$  the eigenvalues of $\Delta_{g,p}$ acting on the subspace $d\Lambda^{p-1}(M)$ of exact $p$-forms, and by $\left\{\lambda''_{k,p}(g); k\ge1 \right\}$ those corresponding to the subspace $d^*\Lambda^{p+1}(M)$ of co-exact $p$-forms, then $\left\{\lambda_{k,p}(g); k\ge1\right\}$ is equal to the re-ordered union of $\left\{\lambda'_{k,p}(g); k\ge1\right\}$ and $\left\{\lambda''_{k,p}(g); k\ge1\right\}$ :
\begin {equation}
\left\{\lambda_{k,p}(g); k\ge1\right\}=\left\{\lambda'_{k,p}(g); k\ge1\right\} \cup \left\{\lambda''_{k,p}(g); k\ge1\right\}.
\end {equation}
Since $d$ and $d^*$ commute with the Laplacian, we have,
$$\lambda'_{k,p}(g)= \lambda''_{k,p-1}(g) $$
and then,
\begin {equation}
\left\{\lambda_{k,p}(g); k\ge1\right\}=\left\{\lambda'_{k,p}(g); k\ge1\right\} \cup \left\{\lambda'_{k,p+1}(g); k\ge1\right\}.
\end {equation}

\medskip

\noindent{\it Proof of Corollary 1.2} : For $p=1$ we have, $\forall k\ge1$,  
$\lambda'_{k,1}(g)=\lambda_{k,0}(g)$ and then  
$$\left\{\lambda_{k,1}(g); k\ge1\right\}=\left\{\lambda_{k,0}(g); k\ge1\right\} \cup \left\{\lambda'_{k,2}(g); k\ge1\right\}.$$
Applying Theorem 1.1 and Korevaar's inequality (3), we get, for any positive integer $N$, the existence of a metric $g_N$ satisfying
$$\lambda'_{1,2}(g_N)> C([g])N^{2\over n}\ge \lambda_{N,0}(g_N).$$
The result of the corollary follows immediately.\hbox{} \hfill $\Box$

\bigskip
\noindent{\it Proof of Theorem 1.1.}
Recall first the idea of the proof in [GP]: they start with the manifold $M$, make a hole and glue 
what they call a cigar, that is the union of 
a standard half-sphere with a cylinder 
$[0,L] \times S^{n-1}$. Since a cigar is diffeomorphic to a ball, the resulting manifold is diffeomorphic to $M$. Then, they endow $M$ with a metric $g_L$ so that it induces the standard metric on the cigar and show that, as $L \to \infty$, $\lambda_{1,p}(g_L)$ is uniformly bounded below by a positive constant for $2 \le p \le n-2$. They 
conclude by noticing that the volume of $(M,g_L)$ goes to infinity with $L$. 

First, let us show that the cigar may be realized as a conformal 
deformation of the $n$-dimensional Euclidean unit ball $B^n$. Indeed, the Euclidean metric $g_{euc}$ has the following expression in polar coordinates 
$$g_{euc}= dr^2 + r^2 d\sigma ^2,$$ 
where $0<r<1$ and $d \sigma ^2$ denotes the standard metric of the 
$(n-1)$-sphere ${\mathbb S}^{n-1}$. For any positive real number $L$, we multiply this metric by the following positive function $f_L$ of the parameter $r$ 

$$f_L(r)= \left \lbrace
\begin{array}{l}

  \frac{4e^{2L}}{(1+e^{2L}r^2)^2} \;\; if \; 0 \le r \le e^{-L}, \\

\\

  \frac{1}{r^2} \;\;\;\;\;\;\;\;\;\;\;\;\; if \; e^{-L} \le r \le 1.

\end{array}
\right.$$
The resulting metric $f_L g_{euc}$ on $B^n$ is so that the part $\left(\{r \le e^{-L} \},f_Lg_{euc}\right) $ is 
isometric to $\left( B^n,\frac{4}{(1+r^2)^2}(dr^2+r^2d\sigma^2)\right)$ that may be identified, via a stereographic projection, with the standard half-sphere (of radius 1), while 
the part $\{e^{-L} < r < 1\}$ endowed with the metric $f_L g_{euc}=(\frac{dr}{r})^2+ 
d\sigma ^2$, is isometric to the standard cylinder $]0,L[ \times S^{n-1}$ with 
the product metric. Hence, the ball $\left(B^n, f_Lg_{euc}\right)$ is isometric to a standard cigar $C_{L}$ of length $L$.

\medskip

Note that the function $f_L$ is only of class $C^1$, 
but, using classical density results, it is possible to smooth it in order to get a smooth metric $\bar f_Lg_{euc}$ conformal to $g_{euc}$ and quasi-isometric to $f_Lg_{euc}$ with a quasi-isometry ratio arbitrarily close to 1.

\medskip

As in [CE], we will need the two following remarks (see the remark O1 and Lemma 2.3 (i) 
p. 343 of [CE]) : 

\begin {itemize} 
\item[(i)] if two Riemannian metrics $g_{1},g_{2}$ on a manifold $M$ are 
$\alpha$-quasi-isometric, then for any smooth positive function $f$, 
the conformal metrics $fg_{1}$ and $fg_{2}$ are also 
$\alpha$-quasi-isometric. 

\smallskip 
\item[(ii)] if $(M,g)$ is a compact Riemannian manifold and $x_{0} \in M$, then 
for any positive $\delta$, there exists a Riemannian metric 
$g_{\delta}$ on $M$ which is flat in a neighborhood of $x_{0}$ and 
$(1+\delta)$-quasi-isometric to $g$ (that is, $\forall v\in TM$, $(1+\delta)^{-2}g_\delta (v,v)\le g(v,v)\le (1+\delta)^2 g_\delta (v,v) $). 
\end {itemize}

Now, let us fix a Riemannian metric $g$ on $M$ and a positive real number $\delta\in \left(0,1\right)$.Thanks to the previous remark (ii),  there exists a metric 
$g_{\delta}$ on $M$ which is flat in a neighborhood of a point $x_{0}\in M$ and 
$(1+\delta)$-quasi-isometric to $g$. Up to a dilatation, we 
can suppose that this neighborhood contains an euclidean ball of radius 
$2$ centered at $x_0$. On the ball of radius 1, $B(x_0,1$), we make the conformal change described above in order
to get a cigar $C_{L}$ of length $L$. On the annulus $\{1\le r \le 2\}$, we interpolate by considering the conformal metric given in polar coordinates by $h(r)(dr^2 + r^2 d\sigma^2)$, where $h$ is any smooth function, $\frac{1}{2} \le h(r) \le 1$, satisfying 

$$h(r)= \left \lbrace
\begin{array}{l}

  \frac{1}{r^2} \;\; if \; 1 \le r \le \frac{4}{3}, \\

\\

  1\;\; if \; \frac{5}{3} \le r \le 2.

\end{array}
\right.$$
Hence, we get, $\forall L>0$, a smooth metric $g_{\delta, L}$ on $M$ conformal to $g_\delta$ and such that a neighborhood of the point $x_0$ is isometric to a standard cigar $C_{L}$ of length $L$. In particular, $(M,g_{\delta, L})$ contains a regular domain $\Omega_L$ naturally identified with the cylinder $[0,L]\times {\mathbb S}^{n-1}$. Moreover, the rest of the manifold $(M-\Omega_L, g_{\delta, L})$ does not depend on $L$ (up to an isometry). Using arguments as in [GP] and [M], it is possible to compare the first eigenvalue $\lambda'_{1,p}(g_{\delta, L})$ of the Laplacian $\Delta_{g_{\delta, L}}$ acting on exact $p$-forms on $M$ with the first positive eigenvalue $\mu_{1,p}(\Omega_L)$ of $\Delta_{g_{\delta, L}}$ acting on exact $p$-forms on $\Omega_L$ under absolute boundary conditions. Indeed, there exist three positive constants $a$, $b$ and $c$, independent of $L$, such that (see [GP], Lemma 1)
$$\lambda'_{1,p}(g_{\delta, L})\ge {a\over {b\over{\mu_{1,p}(\Omega_L)} }+c}. $$
On the other hand, recall that for any integer $p\in [1,n-2]$, the Laplacian acting on $p$-forms on the standard sphere ${\mathbb S}^{n-1}$ admits no zero eigenvalues. Hence, $\forall p\in [2,n-2]$, we have 
\begin{eqnarray}
\nonumber {} \mu_{1,p}(\Omega_L)&=&\mu_{1,p} ([0,L]\times {\mathbb S}^{n-1})\\
 \nonumber{} &\ge& \min \left\{\mu_{1,0}([0,L])+\lambda_{1,p}({\mathbb S}^{n-1}) , \mu_{1,1}([0,L])+\lambda_{1,p-1}({\mathbb S}^{n-1})\right\}\\
 \nonumber{} &\ge& \min \left\{\lambda_{1,p}({\mathbb S}^{n-1}) , \lambda_{1,p-1}({\mathbb S}^{n-1})\right\}\\
\nonumber{} &\ge& \min_{1\le q\le n-2}\lambda_{1,q}({\mathbb S}^{n-1}) .
\end{eqnarray}
Therefore, there exists a positive constant $C$, independent of $L$, such that, $\forall p\in [2,n-2]$ and $\forall L>0$, we have
$$\lambda'_{1,p}(g_{\delta, L})\ge C.$$
The volume one metric $\bar g_{\delta, L}=V(g_{\delta, L})^{-2\over n} g_{\delta, L}$ satisfies, $\forall p\in [2,n-2]$, 
\begin{equation}
\lambda'_{1,p}(\bar g_{\delta, L})\ge C V(g_{\delta, L})^{2\over n} \ge K L^{2\over n}
\end {equation}
for some positive constant $K$. 

Let us denote, $ \forall L>0$, by $\phi_L$ the conformal ratio between $\bar g_{\delta, L}$ and $g_\delta$, that is $\bar g_{\delta, L}=\phi_L g_\delta$. The meric 
 $g_L= \phi_L g$ is then conformal to $g$ and $(1+\delta)$-quasi-isometric to $\bar g_{\delta, L}$ (see remark (i) above). The volume of $g_{_L}$ is bounded below by $(1+\delta)^{-n}V(\bar g_{\delta, L})=(1+\delta)^{-n}$ while its first eigenvalue on exact $p$-forms satisfies $\lambda'_{1,p}(g_{_L})\ge (1+\delta)^{-2(n+2p+1)} \lambda'_{1,p}(\bar g_{\delta, L})$ (see [D] prop. 3.3 p.441). Hence, for any 
$p\in [2,n-2]$, we have 
\begin{equation}
 \lambda'_{1,p}(g_{_L}) V(g_{_L})^{2\over n} \ge \bar K L^{2\over n}
\end {equation}
for some positive constant $\bar K$. 
 
In dimension $n=4$, the space $\Lambda^2(M)$ is invariant by the Hodge star operator $*$ which commutes with the Laplacian and satisfies $ * d\Lambda^1(M) = d^*\Lambda^{3}(M)$. Consequently, using  (5) above, $\lambda_{1,2}(g_{_ L})= \lambda'_{1,2}(g_{_ L}) = \lambda''_{1,2}(g_{_ L})$ and then 
$$\lambda_{1,2}(g_{_ L}) V(g_{_ L})^{2\over n} \ge \bar K L^{2\over n}.$$
In dimension $n\ge 5$, for any $p=2,3,\dots, [{n\over 2}]$,  where $[{n\over 2}]$ is the integer part of $n\over 2$, we have $2\le p<p+1\le n-2$, and then, using (6) and (8) above,  
$$\lambda_{1,p}(g_{_L}) V(g_{_ L})^{2\over n}= \min \left\{\lambda'_{1,p}(g_{_L}), \lambda'_{1,p+1}(g_{_L})\right\} V(g_{_ L})^{2\over n}\ge \bar K L^{2\over n}.$$
This last estimate still holds for $p= [{n\over 2}]+1, \dots, n-2$ since we have $ \lambda_{1,p}(g_{_L})=\lambda_{1,n-p}(g_{_L})$ thanks to the Hodge duality. 

Finally, we have, $\forall n\ge4$ and $\forall p\in [2,n-2]$,
$$\lambda_{1,p}(g_{_L}) V(g_{_L})^{2\over n} \ge \bar K L^{2\over n}$$
and then $\lambda_{1,p}(g_{_L})V(g_{_L})^{2\over n} \to \infty$ as $L\to \infty$. 
Theorem 1.1 follows immediately. \hbox{} \hfill $\Box$

\end{document}